\documentclass{article}
 \usepackage{graphicx, latexsym, epsfig, amssymb}
 \newtheorem{theorem}{Theorem}[section]
\newtheorem{lemma}{Lemma}[section]
 \newtheorem{remark}{Remark}[section]
\newtheorem{example}{Example}[section]
 \newtheorem{definition}{Definition}[section]
\newtheorem{proposition}{Propositin}[section]
 \newtheorem{corollary}{Corollary}[section]
 \def\R{{\mathbb R}}

 \def\Z{{\mathbb Z}}
 \def\d{\displaystyle}

 \title{\bf On the roots of a deformed algebraic system}

 \author{ M.N.~Vrahatis\thanks{vrahatis@math.upatras.gr}\\
 Department of Mathematics,\\
 University of Patras,\\
 GR--26110 Patras, Greece
 \and
 S.~Tanab\'{e}\thanks{tanabe@mccme.ru,
 Partially supported by the Greek State Scholarship Foundation (IKY)}\\
 Department of Mathematics\\
University of Patras,\\
 GR--26110 Patras, Greece\\
 and Department of Mathematics,\\
  Independent University of Moscow,\\
121002 Moscow, Russia}

 \date{}
\begin{document}

%
%
%
%
%
\maketitle
\begin{abstract}
 A problem concerning the shift of roots of a system of algebraic equations is investigated.
 Its conservation and decomposition of a multiple root into simple roots are discussed.


\end{abstract}
\section{Introduction}
The central subject of the present paper is an investigation on the shift
of roots of the system of algebraic equations.
Our central Theorem~\ref{themtineq}
states that the number of real roots of a system located in a compact
set does not change after a sufficiently small perturbation of the system.

As a matter of fact, this kind of fact has been well known to those who study
the deformation of the singularities of  differentiable mappings.
It is, however, a non trivial question how small this perturbation shall
be so that the number of roots in a given compact set remains unchanged.
All the ever existing theorems (see~\cite{ArnoldGV85} \S 12.6)
do not precise the size of the compact set and the perturbation of
the system under question. They state simply that for a compact
set and perturbation,
both of them small enough, the invariance of the number of roots holds.
This situation can be explained by the fact that they
simply treat the notion of local algebra, and consequently
they are valid only in the germ sense.
We try to give here an estimate on the size of admissible
perturbation of the system for a fixed compact set.

Furthermore, we give a result about the decomposition of multiple roots into simple roots.
In particular, our Theorem~\ref{TheoremMult1} assures us the existence of a deformed
system of the original system that possesses only simple roots.

The paper is organized as follows.
In \S 2 preliminary results are presented.
In \S 3 we state that a slightly deformed system has the same number
of zeros as the original system in taking the multiplicities into account.
In \S 4 we give a result about the decomposition of
multiple roots into simple roots.
The paper ends in \S 5 with some concluding remarks.

\section{Preliminary results}
Let us consider the following system of algebraic equations with real
coefficients $a_\alpha^{(i)}$:
\begin{equation}\label{Sys1}
\cases{\begin{array}{l}
f_1(x_1, x_2, \ldots, x_n) = 0, \\[0.2cm]
f_2(x_1, x_2, \ldots, x_n) = 0, \\[0.2cm]
 \hfil\vdots \\[0.2cm]
f_{n}(x_1, x_2, \ldots, x_n) = 0, \\[0.2cm]
\end{array}
}
\end{equation}
with
\begin{equation}
\begin{array}{l}
\d
f_i(x_1, x_2,\ldots,x_n) = \sum_{|\alpha|=m_i} a_\alpha^{(i)} x^\alpha + a_0^{(i)}, \\[0.6cm]
x^\alpha = x_1^{\alpha_1} x_2^{\alpha_2}\cdots x_n^{\alpha_n},\, \kern0.3cm a_0^{(i)}\in \R , \\[0.4cm]
|\alpha| = \alpha_1 + \alpha_2 + \cdots + \alpha_n, \\[0.4cm]
\end{array}
\end{equation}
where the degrees of polynomials are ordered as follows:
\[
m_1 \leqslant m_2 \leqslant \cdots \leqslant m_n.
\]
 Let us consider the situation where the gradient ideal
 $\left\langle\frac{\partial}{\partial x_1}f_{\ell}, \ldots, \right.$
 $\left.\frac{\partial}{\partial x_n}f_{\ell} \right\rangle$
  over $\R[x]$ contains certain power of maximal ideal $\mathfrak{m}^k.$ That is to say,
\begin{equation}\label{eqideal}
\left\langle
\frac{\partial f_{\ell} (x)}{\partial x_1},
\frac{\partial f_{\ell}(x)}{\partial x_2},\ldots,
\frac{\partial f_{\ell}(x)}{\partial x_n}\right\rangle \supset
\left\langle x_1^{\alpha_1}
 x_2^{\alpha_2} \cdots \/ x_n^{\alpha_n}
\right\rangle_{\alpha_1+\cdots+\alpha_n \geqslant k} .
\end{equation}

Let us note by $M_{\alpha}(x) = x_1^{\alpha_1}\cdots x_n^{\alpha_n}$, \
$\alpha_1+\cdots+\alpha_n = k$.
It is possible to consider the set of such monomials as a basis of
$\mathfrak{m}^k.$ The dimension $\mu_n(k)$ of the basis of the ideal
$\mathfrak{m}^k$ can be calculated by the following recurrent relation:
\[
\begin{array}{l}
\displaystyle \mu_n(k) = \sum_{i=0}^k \mu_{n-1}(k-i), \\[0.6cm]
\displaystyle \mu_2(k) = k+1,\\[0.4cm]
\displaystyle \mu_3(k) = \frac{(k+2) (k+1)}{2},\\[0.2cm]
\hfil \vdots
\end{array}
\]
Evidently, $\mu_n(k)$ is the number of the entire lattice points on an
$(n-1)$--dimensional face of the $n$--simplex:
\[
\mu_n(k)= \# \left\{ \left({\alpha_1},\ldots, \alpha_n \right)
\in \Z^n_{\geqslant 0} : \alpha_1+\cdots+\alpha_n =k \right\}.
\]
We recall the notation that:
\[
 \hbox{\em  supp}\,\, \varphi(x) = \{ \alpha \in \Z^n ; \,\,
\varphi_\alpha \neq 0 \},
\]
for a Laurent polynomial $\varphi(x) = \sum_{\alpha \in \Z^n}
\varphi_\alpha x^{\alpha}$.
We call a germ $\varphi(x)$ convenient at zero when the Newton diagram of it at zero contains
non-compact part of all coordinate axes (cf.\ \cite{ArnoldGV85}).
In other words
 $\varphi(x)$ is convenient at zero if it admits the representation,
$$ \varphi(x) = \sum_{i=1}^n
x_i^{\beta_i} + R(x),$$
for $\beta_i \geqslant 1 $ and a certain polynomial $R(x).$
It is easy to see that if $f_{\ell}(x)$ has a convenient germ at zero then their exists
$k {\geqslant} 1$ such that the condition~(\ref{eqideal}) is satisfied.

Suppose that a polynomial vector:
\[
\left(
\matrix{
\varphi_1(x) \cr
\varphi_2(x) \cr
\vdots \cr
\varphi_n(x)\cr
}
\right)
= F \cdot
\left(
\matrix{
\varphi(x) \cr
0\cr
\vdots \cr
0\cr
}
\right),
\]
where $F$ is an invertible constant matrix and
$\varphi(x) \in \mathfrak{m}^{k+1}$ and $\deg \varphi(x) = k' \geqslant k+1$,
i.e.\ $\hbox{\em  supp}\,\, \varphi(x) \subseteq \{
\alpha \in \Z^n_{\geqslant 0}; k+1 \leqslant |\alpha| \leqslant k'\}.$
The question we pose concerns the behavior of roots of a system:
\begin{equation}\label{systemt}
\cases{
\begin{array}{l}
f_1(x) + t \varphi_1(x) = 0, \\[0.2cm]
f_2(x) + t \varphi_2(x) = 0, \\[0.2cm]
\hfil \vdots \\[0.2cm]
f_n(x) + t \varphi_n(x) = 0, \\[0.2cm]
\end{array}
}
\end{equation}
with $t \in [0,1] \subset \R$ as a parameter.

To formulate further statements in a proper way, we introduce notations:
\[
f_i(x) = a_1^{(i)} x^{\vec{v}_1^{(i)}} + a_2^{(i)} x^{\vec{v}_2^{(i)}} + \cdots +
a_{\lambda_i}^{(i)} x^{\vec{v}_{\lambda_i}^{(i)}},
\]
where the vectors
$\vec{v}_j^{(i)} = \left({v}_{j,1}^{(i)},v_{j,2}^{(i)}, \ldots, {v}_{j,n}^{(i)} \right),$
\ $ 1 \leqslant j \leqslant \lambda_i,$ \ satisfy:
\[ \left\langle (1, \ldots, 1),\vec{v}_1^{(i)}\right\rangle =
\left\langle (1, \ldots, 1),\vec{v}_2^{(i)}\right\rangle
= \cdots =
\left\langle (1, \ldots, 1),\vec{v}_{\lambda_i}^{(i)}\right\rangle = m_i .
\]
In general, it is not easy to formulate a sufficient condition on $f_\ell(x)$ so that the
condition~(\ref{eqideal}) holds. We propose here a simple necessary condition for that.

\begin{proposition}
The following isomorphism~(\ref{lattice})
is necessary so that the condition~(\ref{eqideal}) holds,
\begin{equation}\label{lattice}
\begin{array}{r}
\hbox{\rm lattice spanned over $\Z$ of }
\left\{ \vec{v}_1^{(\ell)} - \vec{v}_2^{(\ell)}, \ldots ,
\vec{v}_1^{(\ell)} - \vec{v}_{\lambda_\ell}^{(\ell)} \right\} \cong \\[0.4cm]
\left\{\vec{\alpha} \in \Z^n; \left\langle (1, \ldots, 1), \vec{\alpha}\right\rangle =0 \right\}.
\end{array}
\end{equation}
\end{proposition}

{\bf Proof} If the condition~(\ref{lattice}) does not hold, it is
evidently impossible to create all monomials $x^\alpha$ with
$|\alpha| =k$ as a linear combination of $\frac{\partial}{\partial
x_i}f_\ell(x)$'s. \ {\bf Q.E.D.}

\smallskip

\begin{remark}
This isomorphism can be realized by shifting each lattice point
$\vec{v}$ of the right hand side
of Relation~(\ref{lattice}) towards the lattice point $\vec{v} + (m_i, 0, \ldots, 0)$.
\end{remark}

We give now an example for which condition~(\ref{eqideal}) does not hold.

\begin{example}
Let us consider the following system:
\begin{equation}\label{CounSys}
\cases{
\begin{array}{l}
f_1(x_1, x_2) = a_1^{(1)} x_1^6 + a_2^{(1)} x_1^3 x_2^3 + a_3^{(1)} x_2^6
= 0, \\[0.4cm]
f_2(x_1, x_2) = a_1^{(2)} x_1^{12} + a_2^{(2)} x_1^6 x_2^6 +
a_3^{(2)} x_2^{12} = 0, \\[0.2cm]
\end{array}
}
\end{equation}
where:
\[
\bigl(a_2^{(i)}\bigr)^2 - 4 a_1^{(i)} a_3^{(i)} \neq 0, \ \ i=1,2.
\]
For these polynomials the lattice define on the left hand side of
Relation~(\ref{lattice}) is isomorphic to
\[
\left\{\vec{\alpha} \in \Z^2; \ \ \left\langle (1, 1), \vec{\alpha}\right\rangle =0;
\ \ \frac{\vec{\alpha}}{3} \in \Z^2 \right\} .
\]
Thus, in this case, the condition~(\ref{eqideal}) does not hold.
\end{example}

From now on, we use the notation $f_\ell(x,a)$ instead of $f_\ell(x)$
if we want to emphasize its dependence on the
coefficients $a = (a_1^{\ell},\cdots, a_{\lambda_\ell}^{\ell}).$
For a set of polynomials
${\mathit{\Lambda}}_1^{r} (x, a),
{\mathit{\Lambda}}_2^{r} (x, a), \ldots, {\mathit{\Lambda}}_{\nu_{_{{r}}}}^{r} (x, a) \in \R[x,a]$
homogeneous in variables $a$ we consider the following linear combinations:
\[
{\mathit{\Lambda}}_j^{r+1} (x, a) =
\sum_{i=1}^{\nu_{_{{0}}}} \gamma_i^{(j)} (a)\, x ^{\,\vec{\beta}_i^{\,(j)}} {\mathit{\Lambda}}_i^{r} (x, a),
\ \ \ j=1,2,\ldots {\nu_{_{{r + 1}}}},
\]
where $\vec{\beta}_i^{\,(j)} \in \Z _{\geqslant 0}^n$
and $\gamma_i^{(j)} (a)$ are linear polynomials in variables $a.$

\begin{proposition}\label{chain}
Let as consider the chain of polynomials sets:
\[
\left\{
{\mathit{\Lambda}}_1^{r} (x, a), {\mathit{\Lambda}}_2^{r} (x, a), \ldots,
{\mathit{\Lambda}}_{\nu_{_{r}}}^{r} (x, a)
\right\}, \ \ r =0,  1, 2, \ldots
\]
as above with
\[
{\mathit{\Lambda}}_1^{0} (x, a) =
\frac{\partial}{\partial x_1}f_{\ell}(x,a), \ldots, {\mathit{\Lambda}}_n
^{0} (x, a) =
\frac{\partial}{\partial x_n}f_{\ell} (x,a).
\]
Suppose that for certain $r = L$, some of ${\mathit{\Lambda}}_\ast^{L} (x, a)$'s coincides with
$\lambda^{(s)}(a) M_s (x)$. That is to say there exists
$\bar{h}_1^{(s)} (x, a), \ldots , \bar{h}_n^{(s)} (x, a) \in \R[x, a]$ such that
\[
\lambda^{(s)}(a) M_s (x) = \sum_{i=1}^n \bar{h}_i^{(s)} (x, a) \frac{\partial}{\partial x_i}f_{\ell} (x, a).
\]
Then \ $\deg_a \lambda^{(s)}(a) = L$ \ and \ $\deg_a \bar{h}_i^{(s)} (x, a) = L-1$.
\end{proposition}

{\bf Proof} After the definition of the  recursive process to
create ${\mathit{\Lambda}}_j^{r+1} (x, a) $ from
${\mathit{\Lambda}}_j^{r} (x, a)$ it is clear that
 ${\mathit{\Lambda}}_j^{r} (x, a)$ is a homogeneous polynomial of degree
$r$ in $a.$
The statement is the direct consequence of this fact. \ {\bf Q.E.D.}

\smallskip

\begin{example}\label{example1}
We consider the following example
$$f_\ell(x_1, x_2,a) =f (x,a) = a_1 x_1^5 + a_2 x_1^2 x_2^3 + a_3 x_2^5 = 0,$$
with
$$\frac{\partial}{\partial x_1}f(x,a) = 5a_1 x_1^4 + 2a_2 x_1 x_2^3,$$
$$\frac{\partial}{\partial x_2}f(x,a) = 3a_2 x_1^2 x_2^2+ 5a_3 x_2^4. $$
Then we have the following chain of polynomials
to get $x^{\vec \alpha_2} = x_1^2 x_2^8 $ as a linear combination of
$ \theta_1 = \frac{\partial}{\partial x_1}f(x,a)$  and $\theta_2 = \frac{\partial}{\partial x_2}f(x,a)$:
\[
\begin{array}{l}
\d
 {\mathit{\Lambda}}_1^{1} (x, a) := (2a_2 x_1 x_2^5 -5a_1 x_1^4 x_2^2)\, \theta_1, \\ [0.2cm]
 {\mathit{\Lambda}}_2^{1} (x, a) := (3a_2 x_1^6 -5a_3 x_1^4 x_2^2)\, \theta_2, \\ [0.2cm]
 {\mathit{\Lambda}}_3^{1} (x, a) := 5a_3 x_1^2 x_2^4\, \theta_2, \\[0.4cm]
 {\mathit{\Lambda}}_1^{2} (x, a) := 3a_2 {\mathit{\Lambda}}_1^{1} (x, a),\\ [0.2cm]
 {\mathit{\Lambda}}_2^{2} (x, a) := 3a_2 {\mathit{\Lambda}}_2^{1} (x, a) + 5a_3
                                    {\mathit{\Lambda}}_3^{1} (x, a), \\[0.4cm]
 {\mathit{\Lambda}}_1^{3} (x, a):= 3a_2 {\mathit{\Lambda}}_1^{2} (x, a),\\ [0.2cm]
 {\mathit{\Lambda}}_2^{3} (x, a):= 5a_1 {\mathit{\Lambda}}_2^{2} (x, a),  \\ [0.4cm]
 {\mathit{\Lambda}}_1^{4} (x, a):= 3a_2 {\mathit{\Lambda}}_1^{3} (x, a)+
       5 a_1 {\mathit{\Lambda}}_2^{3} (x, a)= (2^2\;3^3a_2^5 + 5^5a_1^2 a_2^3)x_1^2 x_2^8. \\ [0.2cm]
\end{array}
\]
Thus we have:
\begin{eqnarray*}
& &  \lambda^{(2)}(a) = 2^2\;3^3a_2^5 + 5^5a_1^2 a_2^3, \\
& &  \bar{h}_1^{(2)} (x, a) = (3a_2)^3 (2a_2 x_1x_2^5 -5a_1 x_1^4 x_2^2), \\
& &  \bar{h}_2^{(2)} (x, a) =(5a_1)^2 ((3a_2)^2 x_1^6 -15 a_2 a_3 x_1^4 x_2^2 -(5a_3)^2 x_1^2 x_2^4).
\end{eqnarray*}
For the case of
$x^{\vec \alpha_5} = x_1^5 x_2^5, $ we have:
\begin{eqnarray*}
& &  \lambda^{(5)}(a) = 5a_1(2^23^3a_2^5 + 5^5a_1^2 a_2^3), \\
& & \bar{h}_1^{(5)} (x, a)  =5a_1^2 a_2^3 x_1x_2^5 + 2 \;3^3 \;5 a_1
a_2^4 x_1^4 x_2^2, \\
& &  \bar{h}_2^{(5)} (x, a)  =-2 \cdot 5^2 a_1^2
a_2((3a_2)^2 x_1^6 -15 a_2 a_3 x_1^4 x_2^2 -(5a_3)^2 x_1^2 x_2^4).
\end{eqnarray*}
For the case of $x^{\vec \alpha_{10}} = x_2^{10}, $ we have:
\begin{eqnarray*}
& & \lambda^{(10)}(a) = 5a_3(2^23^3a_2^5 + 5^5a_1^2 a_2^3), \\
& & \bar{h}_1^{(10)} (x, a) =(2^23^3a_2^5 + 5^5a_1^2 a_2^3)x_2^6 -\\
& & \kern2.2cm     -3a_2(5a_1)^2 ( (3a_2)^2 x_1^6 -15 a_2 a_3 x_1^4 x_2^2 -(5a_3)^2 x_1^2 x_2^4), \\
& & \bar{h}_2^{(10)} (x, a)  =-(3a_2)^4 (2a_2 x_1x_2^5 -5a_1 x_1^4 x_2^2).
\end{eqnarray*}
\end{example}

\section{The number of roots of a deformed system}
In this section we state that a slightly deformed system has the same number
of zeros as the original system in taking the multiplicities into account.
We recall here that the index $\ell \in [1, n]$
has been fixed so that $f_{\ell}(x)$ satisfies the condition~(\ref{eqideal}).

\begin{definition}\label{norm}
We introduce the norm:
\[
\| \varphi \| =
\sum_{\alpha \in \, {supp}\, \varphi_{\ell}(x)} | \alpha | \, | \varphi_{\ell, \alpha} |,
\]
where $ \varphi_{\ell} (x) = \sum_{\alpha \in \, {supp}\, \varphi_{\ell}(x)}
\varphi_{\ell, \alpha} x^{\alpha}$.
We name the following value by $ C(a)$:
\[
 C(a) = \max_{1 \leqslant s \leqslant \mu}
 \left(
 \sum_{1 \leqslant j \leqslant n} \, \,
 \sum_{| \vec{\beta}| \leqslant k'-k-1}
 \max_{x \in {\bf K}} \left| h_{j,\ell}^{(s)} (x, a) x^{\vec{\beta}  }\right|
 \right),
\]
for some compact set ${\bf K}$ and $h_{j,\ell}^{(s)} (x, a)= \frac{\bar h_{j}^{(s)} (x, a)}{\lambda^{(s)} (a)}$
after the notation of Proposition~\ref{chain}.\ $\Box$
\end{definition}

\begin{remark}
In general we can not give any reasonable estimate on $C(a)$. In the above Example~\ref{example1}
$h_{j,\ell}^{(s)} (x, a)$ contains coefficients of the form:
\[
\frac{\hbox{polynomial of degree 5 in } (a_1, a_2, a_3)}{a_3 (2^2\, 3^3\, a_2^5 + 5^5\, a_1^2\, a_3^3)}.
\]
This value can be as large as possible if the denominator is very near to zero. The coefficients of
$h_{s,i}^j (x, a)$  contain rational functions in the variable $a$, with denominators $\lambda^{(s)} (a)$
introduced in Proposition~\ref{chain}. \ $\Box$
\end{remark}

Before formulating our main theorem, we recall a simple lemma of linear algebra.

\begin{lemma}\label{invertible}
Let us consider $\mu \times \mu$ matrix $A = \left( a_{ij} \right) \in {\rm End}(\R^{\mu}) $.
If  $ | a_{ij} | < \frac{1}{\mu^2}$, then $({\rm id}_\mu + A)$ is invertible.
\end{lemma}

{\bf Proof} By straightforward calculation of the determinant of
$({\rm id}_\mu + A)$ we have:
\[
\det ({\rm id}_\mu + A) = 1 + a_{11} + a_{22} + \cdots + a_{\mu \mu} + R(a),
\]
where $R(a)$ is a polynomial containing $(\mu^2 - \mu - 1)$ terms of monomials in $(a_{ij})$ whose degrees
are higher than or equal to two and less than or equal to $\mu$.
Evidently, under the condition  $| a_{ij} | < \frac{1}{\mu^2}$ we obtain that $\det({\rm id}_\mu + A) \neq 0$.
\ {\bf Q.E.D.}

\begin{theorem}\label{themtineq}
The number of simple roots of the system~(\ref{systemt}) inside of a compact set
${\bf K}$ coincides with that of the system~(\ref{Sys1})
if $t$ satisfies the following inequality:
\begin{equation}\label{tineq}
t < \frac{1}{ \|\varphi \| \, C(a) \, \mu^2}.
\end{equation}
\end{theorem}

{\bf Proof} Our strategy consists in the construction of a
homotopy that connects the simple  roots of system~(\ref{Sys1})
and those of~(\ref{systemt}).

Suppose that we succeed in constructing a homotopy $x(\tau),$ $0
\leqslant \tau \leqslant t$ with $x(0)=x$ such that
$$ f_s(x(\tau)) +\tau \varphi_s(x(\tau))= f_s(x), \;1\leqslant s
\leqslant n,$$ then the vector field along it satisfies the
following equality:
\[ \frac{d}{d\tau}
\left( \matrix{ f_1(x(\tau)) \cr f_2(x(\tau)) \cr \vdots \cr
f_n(x(\tau))\cr } \right) =\sum_{i=1}^n{\dot x}_i(\tau)
\frac{\partial}{\partial x_i}\left( \matrix{f_1(x(\tau)) \cr
f_2(x(\tau))\cr \vdots \cr f_n(x(\tau))\cr } \right).
\]
In applying this relation to system~(\ref{systemt}), we get,
\[ \sum_{i=1}^n{\dot x}_i(\tau)
\frac{\partial}{\partial x_i}
\left\{
\left(
\matrix{f_1(x(\tau)) \cr
f_2(x(\tau))\hfil\cr \vdots \cr
f_n(x(\tau))\cr } \right) + \tau F \left(
\matrix{ \varphi(x) \cr 0\cr \vdots \cr 0\cr } \right)
\right\}
+ F \left(
\matrix{ \varphi(x) \cr 0\cr \vdots \cr 0\cr } \right)
=  \left(
\matrix{0 \cr 0\cr \vdots \cr 0\cr } \right). 
\]
Further, we shall realize a smooth homotopy
\[
\varphi_{\ell}(x) = \left( \sum_{i=1}^n v_i(x, \tau) \frac{\partial }{\partial x_i}
\right) (f_{\ell}(x) + \tau \varphi_{\ell}(x)).
\]
We remember that we denoted the basis of $\mathfrak{m}^k$ by
$M_{\alpha_i}(x),$ $1 \leqslant i \leqslant \mu = \mu_n(k).$ The
condition~(\ref{eqideal}) entails the following relation:
\begin{equation}\label{Hisum}
\pmatrix{
M_{\alpha_1}(x) \cr
M_{\alpha_2}(x) \cr
 \vdots \cr
M_{\alpha_{\mu}}(x) \cr
} =
H^{(1)} \pmatrix{
\frac{\partial f_1}{\partial x_1} \cr
\frac{\partial f_2}{\partial x_1} \cr
 \vdots \cr
\frac{\partial f_n}{\partial x_1} \cr
}
+
H^{(2)} \pmatrix{
\frac{\partial f_1}{\partial x_2} \cr
\frac{\partial f_2}{\partial x_2} \cr
 \vdots \cr
\frac{\partial f_n}{\partial x_2} \cr
}
+ \cdots +
H^{(n)} \pmatrix{
\frac{\partial f_1}{\partial x_n} \cr
\frac{\partial f_2}{\partial x_n} \cr
 \vdots \cr
\frac{\partial f_n}{\partial x_n} \cr
},
\end{equation}
for some polynomial entry rank--1 $(\mu \times n)$ matrices $H^{(1)}, H^{(2)}, \ldots, H^{(n)}$
of the form:
\[
H^{(i)} = \mu \overbrace{ \cases{
\pmatrix{
0 & \cdots & h_{i,\ell}^{(1)} & 0 & \cdots & 0 \cr
0 & \cdots & h_{i,\ell}^{(2)} & 0 & \cdots & 0 \cr
\vdots & \cdots & \vdots & \vdots & \cdots & \vdots \cr\cr
0 & \cdots & h_{i,\ell}^{(\mu)} & 0 & \cdots & 0 \cr
}}}^n,
\]
where $h_{j,\ell}^{(s)} (x, a)= \frac{\bar h_{j}^{(s)} (x, a)}{\lambda^{(s)} (a)}$
after the notation of Proposition~\ref{chain} concentrated at the $\ell$-th column of the matrix $H^{(i)}$.
One rewrites the relation~(\ref{Hisum}) as follows:
\begin{eqnarray}\label{Hisum2}
\pmatrix{
M_{\alpha_1}(x) \cr
M_{\alpha_2}(x) \cr
 \vdots \cr
M_{\alpha_{\mu}}(x) \cr
} &=&
\left(
\sum_{i=1}^n H^{(i)} \frac{\partial}{\partial x_i}
\right)
\pmatrix{
f_1 + \tau \varphi_1 \cr
f_2 + \tau \varphi_2 \cr
 \vdots \cr
f_n + \tau \varphi_n \cr
} - \nonumber \\[0.2cm]
&& - \tau
\left(
\sum_{i=1}^n H^{(i)} F \frac{\partial}{\partial x_i}
\right)
\pmatrix{
\varphi(x) \cr
0 \cr
 \vdots \cr
0 \cr
}.
\end{eqnarray}
As we supposed that $\varphi(x) \in \mathfrak{m}^{k+1}$, it is easy to see that:
\begin{equation}\label{Hisum3}
\left(
\sum_{i=1}^n H^{(i)} F \frac{\partial}{\partial x_i}
\right)
\pmatrix{
\varphi(x) \cr
0 \cr
 \vdots \cr
0 \cr
} =
A(x)
\pmatrix{
M_{\alpha_1}(x) \cr
M_{\alpha_2}(x) \cr
\vdots \cr
M_{\alpha_{\mu}}(x) \cr
},
\end{equation}
with certain polynomial $(\mu \times \mu)$ matrix $A(x)$, where:
\[
A(x) =
\pmatrix{
g_{1}^{(1)} (x) \, & \, \cdots \, & g_{1}^{(\mu)} (x) \cr
 \vdots & \ddots & \vdots \vspace*{0.2cm}\cr
g_{\mu}^{(1)} (x) \, & \, \cdots \, & g_{\mu}^{(\mu)} (x) \cr
}.
\]
By recalling~(\ref{Hisum2}), we obtain an equation as follows:
\[
(id_\mu + \tau A(x))
\pmatrix{ M_{\alpha_1}(x) \cr M_{\alpha_2}(x) \cr
\vdots \cr
M_{\alpha_{\mu}}(x) \cr }= \left(\sum_{i=1}^n H^{(i)} \frac{\partial}{\partial x_i} \right)
\pmatrix{
f_1 + \tau \varphi_1 \cr f_2 + \tau \varphi_2 \cr
\vdots \cr
f_n + \tau \varphi_n \cr }.
\]
Supposing that $\tau$ is very small, get the inverse to:
\begin{equation}\label{matreq}
(id_\mu + \tau A(x)),
\end{equation}
in the domain $\{ x ; \det (id_\mu + \tau A(x)) \neq 0 \}$.
The inequality~(\ref{tineq}) ensures the invertibility of the
matrix~(\ref{matreq}).
To show this, in view of Lemma~\ref{invertible},
 it is enough to verify that for such a value of
$\tau$ we have:
\begin{equation}\label{taumax}
\tau \left( \max_{1\leqslant i,s \leqslant \mu} \, \, \max_{x \in {\bf K}}
\left| g_s^{(i)} (x) \right| \right)< \frac{1}{\mu^2}.
\end{equation}
In other words, it is enough to prove that:
\begin{equation}\label{taumax2}
\max_{1\leqslant i,s \leqslant \mu} \, \, \max_{x \in {\bf K}}
\left| g_s^{(i)} (x) \right| < \|\varphi \| \, C(a).
\end{equation}
We remember that $\hbox{\em supp}\, g_s^{(i)} (x) \subset
\left\{ \alpha \in \Z^n ; \ k' - m_n \leqslant |\alpha | \leqslant k' - m_i \right\}$.
This is a direct consequence of~(\ref{Hisum3}).
As we have
$\hbox{\em supp}\, \varphi_\ell (x) \subset
\left\{ \alpha \in \Z^n;\right.$ $\left.k + 1 \leqslant \right.$ $\left.|\alpha | \leqslant k' \right\},$
we can find for every $1 \leqslant \lambda \leqslant n$ a series of polynomials
$\xi_{\ell,\lambda}^{(1)} (x), \ldots,$ $\xi_{\ell,\lambda}^{(\mu)} (x)$ such that:
\begin{equation}\label{xi}
\frac{\partial}{\partial x_\lambda} \varphi_\ell(x) =
\sum_{c=1}^{\mu} \xi_{\ell,\lambda}^{(c)} (x) \, M_c (x).
\end{equation}
In terms of these polynomials:
\[
g_s^{(c)} (x) = \sum_{\lambda =1}^n \, h_{\lambda,\ell}^{(s)}
(x) \, \xi_{\ell,\lambda}^{(c)} (x),
\]
and
\[
\hbox{\em supp}\,\xi_{\ell,\lambda}^{(c)} (x) \subset
\left\{ \alpha \in \Z^n ; 0 \leqslant |\alpha | \leqslant k' - k - 1 \right\},
\]
the absolute value of each coefficient of $\xi_{\ell,\lambda}^{(c)} (x)$ can be estimated by
$\| \varphi \| $ after~(\ref{xi}) above. We replace $\xi_{\ell,\lambda}^{(c)} (x)$ by
$\| \varphi \| \times \left( \sum_{\vec{\beta} \in
\hbox{\em supp}\,\xi_{\ell,\lambda}^{(c)} (x)} x^{\vec{\beta}} \right)$
and we get the inequality:
\begin{eqnarray*} 
&
\displaystyle
\kern-3cm
\max_{1\leqslant i,s \leqslant \mu} \, \, \max_{x \in {\bf K}}
\left| g_s^{(i)} (x) \right| <
\left(\sum_{\alpha \in \, {supp}\, \varphi(x)} | \alpha | \, | \varphi_{\ell, \alpha} | \right)\,
\hfill \nonumber \\
&
\displaystyle
\kern-0.3cm
\left( \max_{1\leqslant s \leqslant \mu}
 \sum_{1 \leqslant j \leqslant n} \, \,
 \sum_{ \vec{\beta} \in  \coprod_{i =1}^{n} {supp}\,\frac{\partial \varphi_\ell(x)}{\partial x_i}
 \setminus \coprod_{s =1}^{\mu} {supp}\, M_s(x)}
 \max_{x \in {\bf K}} \left| h_{j,\ell}^{(s)} (x, a)\, x^{\vec{\beta}}\right|
\right),
\end{eqnarray*}
where $A \setminus B = \{ \alpha - \beta \in \Z_{\geqslant 0}^n ; \alpha \in A, \beta \in B \}.$
The Relation~(\ref{xi}) explains the summand of the above inequality.
Therefore, if we set $C(a)$ as in Definition~\ref{norm},
we obtain the inequality~(\ref{taumax2}).
Evidently, $C(a)$ depends not on the coefficients of
$\varphi_1 (x), \ldots, \varphi_n(x)$ but on the powers $k$ and $k'$.
This proves  the invertibility of the matrix~(\ref{matreq}).
Thus,
\begin{equation}\label{relinv}
\pmatrix{ M_{\alpha_1}(x) \cr M_{\alpha_2}(x) \cr
 \vdots \cr
 M_{\alpha_{\mu}}(x) \cr }= (id_\mu + \tau A(x))^{-1}
 \left(\sum_{i=1}^n H^{(i)} \frac{\partial}{\partial x_i} \right)
 \pmatrix{
f_1 + \tau \varphi_1 \cr f_2 + \tau \varphi_2 \cr
\vdots \cr
f_n + \tau \varphi_n \cr }.
\end{equation}
On the other hand:
\[
\left(
\matrix{
\varphi_1(x) \cr
\varphi_2(x) \cr
\vdots \cr
\varphi_n(x)\cr
}
\right)
= F \cdot
\left(
\matrix{
\varphi(x) \cr
0\cr
\vdots \cr
0\cr
}
\right) = F \cdot G
\pmatrix{ M_{\alpha_1}(x) \cr M_{\alpha_2}(x)\cr
\vdots \cr M_{\alpha_{\mu}}(x) \cr }
\]
for some rank--1 $(n\times \mu)$ polynomial matrix:
\[
G=
\pmatrix{
g_1(x) & g_2(x) & \cdots & g_{\mu}(x) \cr
0 & 0 & \cdots & 0 \cr
\vdots & \vdots & \ddots & \vdots \cr
0 & 0 & \cdots & 0 \cr
}.
\]
If we apply $F \cdot G$ from the left to the relation~(\ref{relinv}), we get:
\begin{equation}\label{relinv2}
\pmatrix{
\varphi_1(x) \cr
\varphi_2(x) \cr
\vdots \cr
\varphi_n(x) \cr
} =
F \cdot G (id_\mu + \tau A(x))^{-1} \left(
\sum_{i=1}^n H^{(i)} \frac{\partial}{\partial x_i} \right) \pmatrix{
f_1 + \tau \varphi_1 \cr f_2 + \tau \varphi_2 \cr
\vdots \cr
f_n + \tau \varphi_n \cr },
\end{equation}
or
\begin{equation}\label{relinv3}
\pmatrix{
\varphi(x) \cr
0\cr
\vdots \cr
0 \cr
} =
G (id_\mu + \tau A(x))^{-1} \left(
\sum_{i=1}^n H^{(i)} \frac{\partial}{\partial x_i} \right) \pmatrix{
f_1 + \tau \varphi_1 \cr f_2 + \tau \varphi_2 \cr
\vdots \cr
f_n + \tau \varphi_n \cr }.
\end{equation}
This relation gives rise to an inequality:
\begin{equation}\label{relinv4}
\kern-0.1cm
\pmatrix{
\varphi(x) \cr
0\cr
\vdots \cr
0 \cr
}
\kern-0.1cm=\kern-0.1cm
\sum_{i=1}^{n}\kern-0.1cm
\pmatrix{
0 & \cdots & 0 & {{v_i(x, \tau)}^{^{\kern-0.8cm\vee^{\!\!\!\ell}}}} \kern0.3cm & 0 & \cdots & 0 \cr
0 & \cdots & 0 & 0 & 0 & \cdots & 0 \cr
\vdots & \ddots & \vdots & \vdots & \vdots &  \ddots & 0\cr
0 & \cdots & 0 & 0 & 0 & \cdots & 0 \cr
}
\kern-0.1cm
\frac{\partial }{\partial x_i}
\pmatrix{
f_1 + \tau \varphi_1 \cr f_2 + \tau \varphi_2 \cr
\vdots \cr
f_n + \tau \varphi_n \cr },
\end{equation}
where the matrix in front of the derivative has single non zero $\ell$-th column.
That is to say we obtain the following equalities:
\begin{equation}\label{equalit}
\begin{array}{l}
\displaystyle
\varphi(x) = \left( \sum_{i=1}^n v_i(x, \tau) \frac{\partial }{\partial x_i}
\right) (f_\ell + \tau \varphi_\ell), \\[0.6cm]
\varphi_\ell(x) = F_{\ell 1} \varphi(x)
\end{array}
\end{equation}
The estimate~(\ref{tineq}) ensures that $v_i(x, \tau)$ are real analytic in
$x \in {\bf K}$.
Thus we have constructed a vector field corresponding to the homotopy we need. \ {\bf Q.E.D.}

\smallskip

\begin{example}
Let us consider the following system~\cite{Kearfott79}:
\begin{equation}\label{SystStenger}
\cases{
\begin{array}{l}
f_1(x_1, x_2) = x_1^2 - x_2^2 - 1 = 0, \\[0.4cm]
f_2(x_1, x_2) = x_1^2 + x_2^2 - 2 = 0. \\[0.2cm]
\end{array}
}
\end{equation}
This system has four real solutions
within the square $[-2,2]^2$:
\[
(\pm\sqrt{1.5}, \pm\sqrt{0.5}) \thickapprox (\pm 1.22474487139159, \pm 0.70710678118655).
\]
If we perturb this system with a cubic monomial $\varphi(x)= x_1x_2^2$ as follows:
\begin{equation}\label{SystStenger2}
\cases{
\begin{array}{l}
F_1(x_1, x_2) = x_1^2 - x_2^2 + t x_1x_2^2- 1 = 0, \\[0.4cm]
F_2(x_1, x_2) = x_1^2 + x_2^2 - 2 = 0, \\[0.2cm]
\end{array}
}
\end{equation}
we calculate the constants $C(a)=\frac{5}{2}, \  \| \varphi \| = 3, \mu =2. $
Therefore we have four solutions of the system~(\ref{SystStenger2}) if $t < \frac{1}{2\cdot 3 \cdot 5}.$
In particular, if we use the value $t=0.033<1/30$ by applying the rootfinding method of~\cite{VrahatisI86}
we obtain the following four solutions:
\[
\begin{array}{l}
 (\kern0.3cm1.22054232589618, \kern0.3cm\pm0.71433635683474), \\[0.2cm]
 (-1.22879457180552, \kern0.3cm\pm0.70004564158438).
\end{array}
\]
Furthermore, if we perturb both equations of the system~(\ref{SystStenger}) with the same cubic monomial
$\varphi(x)= x_1x_2^2$ as follows:
\begin{equation}\label{SystStenger3}
\cases{
\begin{array}{l}
F_1(x_1, x_2) = x_1^2 - x_2^2 + t x_1x_2^2 - 1 = 0, \\[0.4cm]
F_2(x_1, x_2) = x_1^2 + x_2^2 + t x_1x_2^2 - 2 = 0, \\[0.2cm]
\end{array}
}
\end{equation}
we again calculate the constants $C(a)=\frac{5}{2}, \  \| \varphi \| = 3, \mu =2. $
Therefore we have four solutions of the system~(\ref{SystStenger3}) if $t < \frac{1}{2\cdot 3 \cdot 5}.$
In particular, if we use the value $t=0.033<1/30$ by applying the rootfinding method of~\cite{VrahatisI86}
we obtain the following four solutions:
\[
\begin{array}{l}
 (\kern0.3cm1.21652265747566, \kern0.3cm\pm0.70710678118655), \\[0.2cm]
 (-1.23302265747566, \kern0.3cm\pm0.70710678118655).
\end{array}
\]
\end{example}

\smallskip

\begin{example}
Let us consider the following system:
\begin{equation}\label{Sys3d}
\cases{
\begin{array}{l}
f_1(x_1, x_2, x_3) = \left(x_1^2 + \d \frac{\d x_2^2}{\d 4} - x_3^2\right)
\left( \d \frac{x_1^2}{4} + x_2^2 - x_3^2\right) - \d \frac{\d x_3^4}{\d 81} = 0, \\[0.4cm]
f_2(x_1, x_2, x_3) = (x_1 + x_2)^2 + 36 (x_1 - x_2)^2 - 9 x_3^2 = 0, \\[0.4cm]
f_3(x_1, x_2, x_3) = \d \frac{\d x_1^2}{\d 4} + x_2^2 + \d \frac{\d x_3^2}{\d 9} -1 =0.
\end{array}
}
\end{equation}
By applying the rootfinding method of~\cite{VrahatisI86}
we obtain the following sixteen real solutions
within the cube $[-2,2]^3$:
\[
\begin{array}{l}
 (\kern0.3cm0.62830967308983, \kern0.3cm0.91412675198426, \kern0.1cm\pm0.76883755100759), \\
 (-0.62830967308983, -0.91412675198426, \kern0.1cm\pm0.76883755100759), \\
 (\kern0.3cm0.49635596537865, \kern0.3cm0.91441703848857, \kern0.1cm\pm0.95929271740718),\\
 (-0.49635596537865, -0.91441703848857, \kern0.1cm\pm0.95929271740718), \\
 (\kern0.3cm1.11731818404380, \kern0.3cm0.76796989195429, \kern0.1cm\pm0.93973420474984),\\
 (-1.11731818404380, -0.76796989195429, \kern0.1cm\pm0.93973420474984),\\
 (\kern0.3cm1.22450432822695, \kern0.3cm0.66467487192937, \kern0.1cm\pm1.28459776563576),\\
 (-1.22450432822695, -0.66467487192937, \kern0.1cm\pm1.28459776563576).
\end{array}
\]
We observe that these roots are invariant under the actions  of  a group
$G:=(\Z/2\Z) \times (\Z/2\Z)$ generated by two generators $(x_1,x_2,x_3)\mapsto(-x_1,-x_2,x_3)$
and $(x_1,x_2,x_3)\mapsto(x_1,x_2,-x_3)$ due to the invariance of the system~(\ref{Sys3d})
itself under the same group action.
If we perturb this system with a quadratic monomial $\varphi(x)= x_2^2$ as follows:
\begin{equation}\label{Sys3d2}
\cases{
\begin{array}{l}
f_1(x_1, x_2, x_3) = \left(x_1^2 + \d \frac{\d x_2^2}{\d 4} - x_3^2\right)
\left( \d \frac{x_1^2}{4} + x_2^2 - x_3^2\right) - \d \frac{\d x_3^4}{\d 81} = 0, \\[0.4cm]
f_2(x_1, x_2, x_3) = (x_1 + x_2)^2 + 36 (x_1 - x_2)^2 - 9 x_3^2  + t x_2^2 = 0, \\[0.4cm]
f_3(x_1, x_2, x_3) = \d \frac{\d x_1^2}{\d 4} + x_2^2 + \d \frac{\d x_3^2}{\d 9} -1 =0,
\end{array}
}
\end{equation}
we calculate the constants $C(a)=\frac{1}{2}, \  \| \varphi \| = 2, \mu =3. $
Therefore we have sixteen solutions of the system~(\ref{Sys3d2}) if $t < \frac{1}{3^2}.$
In particular, if we use the value $t=0.1<1/9$ by applying the rootfinding method of~\cite{VrahatisI86}
we obtain the following  sixteen real solutions within the cube $[-2,2]^3$:
\[
\begin{array}{l}
 (\kern0.3cm0.63087661393950, \kern0.3cm0.91351892559324, \kern0.1cm\pm0.77060795720733), \\
 (-0.63087661393950, -0.91351892559324, \kern0.1cm\pm0.77060795720733), \\
 (\kern0.3cm0.49896002229193, \kern0.3cm0.91405620623649, \kern0.1cm\pm0.95934810309529),\\
 (-0.49896002229193, -0.91405620623649, \kern0.1cm\pm0.95934810309529),\\
 (\kern0.3cm1.11568183127565, \kern0.3cm0.76857206607484, \kern0.1cm\pm0.93967783245553),\\
 (-1.11568183127565, -0.76857206607484, \kern0.1cm\pm0.93967783245553),\\
 (\kern0.3cm1.22357424633595, \kern0.3cm0.66527556809517, \kern0.1cm\pm1.28379297777855),\\
 (-1.22357424633595, -0.66527556809517, \kern0.1cm\pm1.28379297777855).
 \end{array}
\]
One remarks here also the invariance of the roots
 under the above mentioned group $G$ due to the invariance of the system~(\ref{Sys3d2}) itself.\ {$\Box$}
\end{example}

\smallskip

As for the equation~(\ref{Sys1}) we establish the following
theorem:

\begin{theorem}\label{thmKov}
Let us consider a system of algebraic equations obtained as a perturbation of~(\ref{Sys1}):
\begin{equation}\label{Sys2}
\cases{
\begin{array}{l}
F_1(x_1, x_2, \ldots, x_n) = 0, \\[0.2cm]
F_2(x_1, x_2, \ldots, x_n) = 0, \\[0.2cm]
  \hfil\vdots \\[0.2cm]
F_{n}(x_1, x_2, \ldots, x_n) = 0. \\
\end{array}
}
\end{equation}
Suppose that on a ball\/ ${\bf B}_{r} = \{x \in \R^n; \, \, |x| \leqslant r \}$
we have:
\begin{equation}\label{rank}
{\rm rank} \left( \frac{\partial}{\partial x_j} f_k (x) \right)_{1\leqslant j,k \leqslant n} = n.
\end{equation}
Furthermore, we impose a condition on $(F_1, F_2, \ldots, F_{n})$:
\begin{equation}\label{Cond}
\max_{x,y \in {\bf B}_{r}}\kern-0.1cm
\left| \hbox{vector component of }\kern-0.1cm
\left( \frac{\partial}{\partial y_j} f_k (y) \right)_{\kern-0.1cm1\leqslant j,k
\leqslant n}^{-1} \kern-0.1cm
\left[
\begin{array}{l}
(F_1 - f_1) (x) \\[0.2cm]
(F_2 - f_2) (x) \\[0.2cm]
\hfil \vdots\\[0.2cm]
(F_{n} - f_{n}) (x)
\end{array}
\right]
\right| < \varepsilon,
\end{equation}

\begin{equation}\label{CondF}
\max_{x,y \in {\bf B}_{r}}\kern-0.1cm
\left| \hbox{vector component of }\kern-0.1cm
 \left( \frac{\partial}{\partial y_j} F_k (y) \right)_{\kern-0.1cm1\leqslant j,k
\leqslant n}^{-1} \kern-0.1cm
\left[
\begin{array}{l}
(F_1 - f_1) (x) \\[0.2cm]
(F_2 - f_2) (x) \\[0.2cm]
\hfil \vdots\\[0.2cm]
(F_{n} - f_{n}) (x)
\end{array}
\right]
\right| < \varepsilon,
\end{equation}
where $\varepsilon$ is strictly less than distance of any root of~(\ref{Sys1}) in ${\bf B}_{r}$
to the boundary $\partial\,{\bf B}_{r}$. Suppose that the system~(\ref{Sys1}) has no multiple real
roots.

Under these assumption the equality:
\[
\#\{\hbox{\textmd{\kern-0.05cmreal simple roots of }} (\ref{Sys1})\kern-0.05cm \hbox{ in }\kern-0.05cm {\bf B}_{r}\}
\kern-0.07cm=\kern-0.07cm
\#\{\hbox{\textmd{\kern-0.05cmreal simple roots of }} (\ref{Sys2})\kern-0.05cm \hbox{ in }\kern-0.05cm {\bf B}_{r}\}
\]
holds.
\end{theorem}

{\bf Proof} We solve the homotopy equation with respect to smooth
diffeomorphism
\[
  x_i + h_i(x,\tau), \kern0.5cm 0 \leqslant i\leqslant n,
\]
that satisfies
\begin{eqnarray}
&& h_i(x,0)  =  0, \\
&& f_k (x_0+h_0(x,\tau),\ldots,x_n+h_n(x,\tau))  = f_k(x_0,\ldots,x_n) + \nonumber \\
&&\kern4cm +\tau (F_k(x)-f_k(x)), \ \ \ 0\leqslant \tau \leqslant 1.
\end{eqnarray}
The system gives rise to a system of $(n+1)$ nonlinear differential
equations:
\begin{equation}\label{diffeq}
\sum_{j=0}^n \frac{\partial h_j}{\partial \tau}
\partial_j f_k(x_0+h_0(x,\tau),\ldots,x_n+h_n(x,\tau))=F_k(x)-f_k(x),
\end{equation}
for $k=0,1,\ldots,n.$
From the assumption~(\ref{rank}), Eq.~(\ref{diffeq}) is always
solvable in the class of real analytic functions  so far
as
\[
\det\left( \left( \frac{\partial}{\partial x_j} f_k (x)
\right)_{1\leqslant j,k \leqslant n}\right) \neq 0
\]
after Cauchy--Kovalevskaya's theorem~\cite{Kovalevskaya875} on the quasi-linear partial differential equation.
After the conditions~(\ref{Cond})-(\ref{CondF}) and Eq.~(\ref{diffeq}), $| \frac{\partial h_j}{\partial \tau}|$
is always strictly less than $\varepsilon$. Therefore $| h_j(\tau)| < \varepsilon\, \tau$ and
$|x+ h_j(\tau)| < |x|+ \varepsilon \, \tau < r$
for $x$ root of~(\ref{Sys1}) located in the ball ${\bf B}_r$. Thus the homotopy equation admits a real
analytic solution that connects
$x \in {\bf B}_{r}$ with  $x+h(x,1) \in {\bf B}_{r}$. \ {\bf Q.E.D.}

\begin{corollary}\label{Corollary}
If the conditions on the analyticity of the homotopy constructed in the above Theorems~\ref{themtineq}
and~\ref{thmKov} are fulfilled then none of the roots of the deformed system~(\ref{defsys})
crashes with another and consequently  no new multiple roots are
created after the proposed perturbation.
\end{corollary}

{\bf Proof} Assume that after the proposed perturbation a multiple
root is created. Then the homotopy constructed in the above
Theorems~\ref{themtineq} and~\ref{thmKov} looses its analyticity
with respect to the parameter~$\tau$. \ {\bf Q.E.D.} \noindent

\section{Decomposition of multiple roots}
In this section we recall facts about the decomposition of
multiple roots into simple roots.

\begin{definition}
For the system~(\ref{Sys1}), we define the Jacobian function:
\begin{equation}
(jf)(x) = \det
\left[
\matrix{
\d\frac{\partial f_1}{\partial{x_1}} \, \,    \cdots \, \,  \d \frac{\partial
f_1}{\partial{x_n}}\cr
\vdots \hfil\ddots \hfil \vdots \cr\cr
\d \frac{\partial f_{n}}{\partial{x_1}}  \, \,   \cdots \, \,  \d \frac{\partial
f_{n}}{\partial{x_n}} \cr}
\right].
\end{equation}
Let us denote by $V_{jf} = \{ x \in \R^{n} : jf(x)=0 \}$ the zero
set of $jf(x)$.
We use the notation $Q_{jf}$ for the following set \[
  V_{jf} \cap \{x : (f_1^2 + \cdots + f_n^2 )(x) =0 \}.
\]
\end{definition}

Then we have the following result:

\begin{theorem}
If $ Q_{jf}= \emptyset$
then the system~(\ref{Sys1}) has only simple roots, while if
$ Q_{jf} \neq \emptyset$
then the system~(\ref{Sys1}) has multiple roots.
\end{theorem}

{\bf Proof} The proof follows from well known facts in the
singularity theory~\cite{ArnoldGV85}. We prove the contrapositive
of the statements. Namely we can easily see  that the existence of
multiple roots yields non-emptiness of  $Q_{jf}$. On the other
hand, the existence of only simple roots
 entails the emptiness of it $Q_{jf}$.
\ {\bf Q.E.D.}

\begin{definition}
Let us
denote by $(Jf)$ a vector valued ideal
\begin{equation}
(Jf) = \left\langle
\left[
\matrix{
\d\frac{\partial f_1}{\partial{x_1}} \, \,    \cdots \, \,  \d \frac{\partial
f_1}{\partial{x_n}}\cr
\vdots \hfil\ddots \hfil \vdots \cr\cr
\d \frac{\partial f_{n}}{\partial{x_1}}  \, \,   \cdots \, \,  \d \frac{\partial
f_{n}}{\partial{x_n}} \cr}
\right]
\left[
\matrix{
\R[x] \cr
\vdots \cr
\R[x] \cr}
\right]
 \right\rangle.
\end{equation}
\end{definition}

\begin{theorem}\label{TheoremMult1}
Let us consider a system like~(\ref{Sys1}) for which we know that it possesses
$m$ real roots
with multiplicities $n_1, n_2,\ldots, n_m$ where $1\leqslant  n_j$ for
$1\leqslant j \leqslant m$.
Then there exists a vector polynomial:
\begin{equation}\label{vectpol}
\left[ \left.
\begin{array}{l}
H_1(x) \\
\hfil \vdots \\
H_{n}(x)
\end{array}
\right] \in (\R [x])^{n} \right/(Jf)
\end{equation}
such that the system of equations
\begin{equation}\label{defsys}
\begin{array}{l}
(f_1 + H_1)(x) = 0, \\
\hfil \vdots \\
(f_{n} + H_{n})(x) = 0,
\end{array}
\end{equation}
has $n_1 + n_2 + \cdots + n_m$ simple real roots.
\end{theorem}

{\bf Proof} We remark that
$$(f_i + H_i)(x) = \sum_{k=0}^{m_i} h_{i,k}(x')
x_1^{m_i-k},$$
with $h_{i,m_i}(x') \not \equiv 0, $$h_{i,0}(x') \not \equiv 0,$
$x'=(x_2, \cdots, x_n),$
after certain permutation of variables $x.$
It is well known that, there exists a perturbation $H_i(x)$
such that the equation $\sum_{k=0}^{m_i} h_{i,k}(x')
x_1^{m_i-k}=0$ has $m_i$ simple roots for a codimension $1$
set of $x'$ (cf.\ \cite{ArnoldGV85}).This fact entails that the system~(\ref{defsys})
also possesses as much simple roots as~(\ref{Sys1}) has. \ {\bf Q.E.D.}

\begin{remark}
One can understand this theorem by an intuitive way.
Let us denote by $I_i :={supp (f_i + H_i)}$
the set of powers present in the polynomial $f_i + H_i.$
If the discriminant of the system
\begin{equation}\label{defsys2}
(f_i + H_i)(x) = \sum_{\alpha \in I_i } f_{\alpha}
x^{\alpha},   \; 1 \leqslant i \leqslant n,
\end{equation}
say, $\Delta(f_\alpha) \in {\R}[f_{\alpha_1}, \ldots, f_{\alpha_{\sum_{i=1}^n|I_i|}}]$
does not vanish, then the roots of the system~(\ref{defsys2})
are all simple. That is to say the set of the coefficients
of the system~(\ref{defsys2}) for which the system has multiple roots
is of codimension one in the space of coefficients ${\R}^{\sum_{i=1}^n|I_i|}.$
This fact is known under the name of Bertini-Sard theorem~\cite{ArnoldGV85}.
\end{remark}

\begin{remark}
It is worthy to notice that one shall choose a proper vector polynomial
(\ref{vectpol}),
to get distinct simple roots for the deformed system (\ref{defsys}). For example
if $f_1 = x^3$ the $f_1=0$ has a triple root at $x=0$. If we take $H_1 = -x$
then
$f_1+H_1 = x^3-x = x(x^2-1)=0$ has 3 distinct roots $x=-1,0,1$, while for
$H_1=x$, the
equation $x^3 +x =0$ has 1 real simple root at $x=0$ and two distinct complex
roots.
\end{remark}

\begin{example}
Let us consider the following system:
\begin{equation}\label{SystStengerMult}
\cases{
\begin{array}{l}
f_1(x_1, x_2) = x_1^2 - x_2^2 - 1 = 0, \\[0.4cm]
f_2(x_1, x_2) = x_1^4 + x_2^2 - 1 = 0. \\[0.2cm]
\end{array}
}
\end{equation}
This system has two multiple real solutions $(\pm1,0)$
within the square $[-2,2]^2$.

If we perturb this system with $H_1=0$ and the simple linear polynomial $H_2 = t (x_1 - 2)$ where $0 < t \leqslant 0.5$
as follows:
\begin{equation}\label{SystStengerMult2}
\cases{
\begin{array}{l}
(f_1+H_1)(x_1, x_2) = x_1^2 - x_2^2 - 1 = 0, \\[0.4cm]
(f_2+H_2)(x_1, x_2) = x_1^4 + x_2^2 + t (x_1 - 2) - 1 = 0, \\[0.2cm]
\end{array}
}
\end{equation}
then we have four simple real solutions.
In particular, if we use the value $t=0.5$ by applying the rootfinding method of~\cite{VrahatisI86}
we obtain the following four real simple solutions:
\[
  (1.07123233675477, \pm0.38410769233261),
\]
\[
  (-1.20970135357686, \pm0.68071827127359).
\]
While if we use the value $t=0.025$ we obtain the following four real simple solutions:
\[
  (1.00412951827050,  \pm0.09097301502177),
\]
\[
  (-1.01237171332486, \pm0.15778620326351).
\]
Finally, if we use the value $t=0.0125$ we obtain the following four real simple solutions:
\[
  (1.00207398824224,  \pm0.06443817123186),
 \]
\[
  (-1.00621769007449, \pm0.11168724107454).
\]
\end{example}

\begin{example}
Let us consider the following system:
\begin{equation}\label{Syst3DMult}
\cases{
\begin{array}{l}
f_1(x_1, x_2, x_3) = \left((x_1-x_2)^3-x_3^2 (x_1+x_2-x_3)\right)\\[0.2cm]
\kern2.7cm     \left((-x_1-x_2)^3-x_3^2 (-x_1+x_2-x_3)\right) \\[0.2cm]
\kern2.7cm     \left((-x_1+x_2)^3-x_3^2 (-x_1-x_2-x_3)\right) \\[0.2cm]
\kern2.7cm     \left((x_1+x_2)^3-x_3^2 (x_1-x_2-x_3)\right) = 0, \\[0.4cm]
f_2(x_1, x_2, x_3) = x_1^2+x_2^2- \d \frac{\d x_3^2}{\d 2} = 0, \\[0.4cm]
f_3(x_1, x_2, x_3) = x_1^2 + \d \frac{\d x_2^2}{\d 9} + \d \frac{\d x_3^2}{\d 4} -1 = 0. \\[0.2cm]
\end{array}
}
\end{equation}
This system has eight simple real solutions:
\[
 (\pm0.25926718242254, \pm1.21300057180546, \pm1.75418919109753),
\]
and eight triple real solutions:
\[
 (\pm0.68824720161168, \pm0.68824720161168, \pm1.37649440322337).
\]
within the cube $[-2,2]^3$.
We observe that these roots are invariant under the actions  of  a group
$\Gamma:=(\Z/2\Z) \times (\Z/2\Z)\times (\Z/2\Z)$ generated by three generators $(x_1,x_2,x_3)\mapsto(-x_1,x_2,x_3),$
$(x_1,x_2,x_3)\mapsto(x_1,-x_2,x_3)$ and $(x_1,x_2,x_3)\mapsto(x_1,x_2,-x_3)$ due to the invariance of
the system~(\ref{Syst3DMult})
itself under the same group action.

If we perturb this system with $H_1$ where $0 < t \leqslant 0.5$
and $H_2 =H_3= 0$ as follows:
\begin{equation}\label{Syst3DMultDIAT}
\cases{
\begin{array}{l}
(f_1 +H_1)(x_1, x_2, x_3) = \\[0.2cm]
\kern1.cm     \left((x_1-x_2)^3 -t(x_1-x_2)-x_3^2 (x_1+x_2-x_3)\right)\\[0.2cm]
\kern1.cm     \left((-x_1-x_2)^3-t(-x_1-x_2) -x_3^2 (-x_1+x_2-x_3)\right) \\[0.2cm]
\kern1.cm     \left((-x_1+x_2)^3-t(-x_1+x_2) -x_3^2 (-x_1-x_2-x_3)\right) \\[0.2cm]
\kern1.cm     \left((x_1+x_2)^3-t(x_1+x_2)-x_3^2 (x_1-x_2-x_3)\right) = 0, \\[0.4cm]
(f_2+H_2)(x_1, x_2, x_3) = x_1^2+x_2^2- \d \frac{\d x_3^2}{\d 2} = 0, \\[0.4cm]
(f_3+H_3)(x_1, x_2, x_3) = x_1^2 + \d \frac{\d x_2^2}{\d 9} + \d \frac{\d x_3^2}{\d 4} -1 = 0, \\[0.2cm]
\end{array}
}
\end{equation}
then we have thirty two simple real solutions.
In particular, if we use the value $t=0.5$
we obtain the following eight real solutions which are shifts of the simple solutions to System~(\ref{Syst3DMult}):
\[
  (\pm0.27142016486929, \pm1.20645760731621, \pm1.74883324771051).
\]
Also we obtain the following twenty four simple real solution:
\[
 (\pm0.68824720161168, \pm0.68824720161168, \pm1.37649440322337),
\]
\[
 (\pm0.78897550317143, \pm0.32932116069209, \pm1.20907797224513),
\]
\[
 ( \pm0.44474589932680, \pm1.07278013064881, \pm1.64234961179579),
\]
that originate from the triple solutions to System~(\ref{Syst3DMult}).
We remark that the first ones of the above solutions coincide with the triple solutions
to the original system. These roots are also invariant under actions of the group $\Gamma.$
\ {$\Box$}
\end{example}

\section{Concluding remarks}

A problem concerning the shift of roots of a system of algebraic equations has been investigated.
Its conservation and decomposition of a multiple root into simple roots have been discussed.

To this end, with our central Theorem~\ref{themtineq} we
show that the number of real roots of a system located in a compact
set does not change after a sufficiently small perturbation of the system.
This theorem can be applied to high dimensional CAD where it is
sometimes needed to calculate intersection of several hypersurfaces
that are perturbation of a set of original (unperturbed) hypersurfaces.
For example, to draw a 3D (three dimensional) picture
of a real algebraic surface that obtained  as a deformation of a known one,
the question of the shift of roots plays essential role~\cite[\S 9.6]{PatrikalakisM02}.
We hope that our results to this direction represent certain interests to those
who are interested in the application of algebraic equations to
computer graphics.

Furthermore, we give a result about the decomposition of multiple roots into simple roots.
In particular, our Theorem~\ref{TheoremMult1} assures us the existence of a deformed
system~(\ref{defsys}) of the original system~(\ref{Sys1}) that possesses only
simple roots.
This result can be used in many cases including the computation of
the topological degree~\cite{Kearfott79,MourrainVY02,Stenger75,Stynes79,Stynes81}
in order to examine the solution set of a system of equations and to
obtain information on the existence of solutions, their number and their
nature~\cite{AlexandroffH65,KavvadiasV96,MourrainVY02,Picard892,Picard22}.

%
%

\end{document}